\documentclass[12pt,A4]{amsart}
\usepackage{amsmath}
\usepackage{amsthm}
\usepackage{tikz}
\usepackage{tikz-cd}
\usepackage{amssymb}
\usepackage{graphicx}
\usepackage{thmtools}
\usepackage{thm-restate}
\usepackage{setspace}
\usepackage{gensymb}
\usepackage{cleveref}
\usepackage{fullpage}
\usetikzlibrary{patterns}
\newtheorem{theorem}{Theorem}[section]

\theoremstyle{definition}
\newtheorem{definition}{Definition}[section]
\newtheorem{example}{Example}[section]

\allowdisplaybreaks
\usepackage{comment}
\usepackage{xcolor}

\allowdisplaybreaks[1]

\begin{document}
\title{Identifying logarithmic tracts}
\author{James Waterman}
\address{School of Mathematics and Statistics, The Open University, Walton Hall, Milton Keynes MK7 6AA, UK}
\curraddr{}
\email{james.waterman@open.ac.uk}
\date{}
\maketitle 
\begin{abstract}
We show that a direct tract bounded by a simple curve is a logarithmic tract and further give sufficient conditions for a direct tract to contain logarithmic tracts. As an application of these results, an example of a function with infinitely many direct singularities, but no logarithmic singularity over any finite value, is shown to be in the Eremenko-Lyubich class.
\end{abstract}
\section{Introduction}

Let $D$ be an unbounded domain in $\mathbb{C}$ whose boundary consists of piecewise smooth curves and suppose that the complement of $D$ is unbounded. Further, let $f$ be a complex valued function whose domain of definition contains the closure $\overline{D}$ of $D$. Then $D$ is called a  \textit{direct tract} of $f$ if $f$ is holomorphic in $D$, continuous in $\overline{D}$, and if there exists $R>0$ such that $|f(z)|=R$ for $z \in \partial D$ while $|f(z)| > R$ for $ z \in D$. We call $R$ the \textit{boundary value} of the direct tract. If, in addition, the restriction $f:D \rightarrow \{z \in \mathbb{C}: |z| > R\}$ is a universal covering, then $D$ is called a \textit{logarithmic tract}. Further, if $f'(z)=0$ we say $z$ is a \textit{critical point} and that $f(z)$ is a \textit{critical value}. Finally, if $\Gamma:[0,\infty) \rightarrow \mathbb{C}$ with $|\Gamma(t)| \rightarrow \infty$ as $t \rightarrow \infty$ and $f(\Gamma(t))\rightarrow a$ as $t\rightarrow \infty$, then $a \in \hat{\mathbb{C}}$ is an \textit{asymptotic value} of $f$ associated with the \textit{asymptotic path} $\Gamma$.

Every transcendental entire function has a direct tract. Moreover, all direct tracts of a function in the Eremenko-Lyubich class $\mathcal{B}$ are logarithmic tracts, for sufficiently large $R$ in the above, where the class $\mathcal{B}$ is the class of transcendental entire functions for which the set of critical and asymptotic values is bounded. 

Logarithmic tracts are useful in proving a variety of results.  Bara\'nski, Karpi\'nska, and Zdunik     \cite{BKZ09} showed that, if a meromorphic function $f$ has a logarithmic tract, then the dimension of the Julia set of $f$ is strictly greater than $1$.  Rottenfusser, R\"{u}ckert, Rempe,
              and Schleicher \cite{RRRS}, as well as Bergweiler, Rippon, and Stallard \cite{Tracts}, proved results on the structure of the escaping set for functions with a logarithmic tract. Also, one can construct orbits of points that escape slower than any given sequence within a logarithmic tract \cite{Waterman19}. These results use the fact that within a logarithmic tract one can define a logarithmic transform of $f$, which gives an expansion estimate due to Eremenko and Lyubich \cite[Lemma~1]{EL92}. 
Therefore, it is useful to be able to identify such logarithmic tracts. In this paper, we first give a geometric condition for a direct tract to be logarithmic.

\begin{theorem}\label{Single Curve}
Let $D$ be a direct tract whose boundary is an unbounded simple curve. Then $D$ is logarithmic.
\end{theorem}

While direct tracts bounded by a single curve are logarithmic, logarithmic tracts need not be bounded by a single curve. The function $e^{e^z}$ has direct tracts which are horizontal strips when the boundary value of the direct tract considered is $R=1$. No critical points lie in these direct tracts and the asymptotic values do not lie in the image of the tracts, so these direct tracts are logarithmic. However, with the additional assumption that there are no asymptotic paths in a logarithmic tract $D$ with asymptotic values having the same modulus as the boundary value of the direct tract, then $D$ will be bounded by a single curve. Further, if $D$ is a logarithmic tract for $|f(z)|=R$ on $\partial D$, then for all $R'>R$, $D$  will be bounded by a single curve.

In the case where a direct tract is bounded by more than one curve, and possibly by infinitely many curves, we give sufficient conditions for the direct tract to \textit{contain} logarithmic tracts. If a simply connected direct tract is bounded by finitely many curves, these results imply that the direct tract contains only logarithmic tracts and asymptotic paths with asymptotic values of modulus equal to the boundary value of the direct tract. Further,  in the case where there are finitely many boundary curves of a simply connected direct tract, there can be only finitely many critical points in the direct tract, at most $m-1$, where $m$ is the number of logarithmic tracts contained in it.
In order to better describe this situation, we define an access to a point, as is done in \cite{Accesses}, and a channel of a tract, a new concept. Note that, only the concept of an access to \textit{infinity} will be used.
\begin{definition}
Let $U$ be a simply connected domain in $\mathbb{C}$. Fix a point $z_0 \in U$. A homotopy class of curves $\gamma : [0, 1] \rightarrow \hat{\mathbb{C}}$ such that $\gamma([0,1)) \subset U$, $\gamma(0)=z_0$, and $\gamma(t)\rightarrow \infty$ as $t\rightarrow 1$ is called an \textit{access} from $U$ to $\infty$.
\end{definition}

\begin{definition}
Let $D$ be a direct tract. An unbounded simply connected component $G$ of $ \{z \in D: |z|>r\}$ for $r>0$  is called a \textit{channel} of $D$ if there exists exactly one access to $\infty$ in~$G$.
\end{definition}

Note that such a channel must be bounded by a single unbounded simple curve.

\begin{theorem}\label{InftyCurves}
Let $D$ be a direct tract  of $f$, bounded by more than one, and possibly infinitely many, distinct unbounded simple curves, with $|f(z)|=R$, for $z\in \partial D$. Then, for any channel in $D$, either
\begin{itemize}
\item  the channel contains a logarithmic tract, or 
\item $f$ has the same finite asymptotic value of modulus $R$ along all paths to infinity in the channel.
\end{itemize}

In particular, if $D$ is bounded by $m$  distinct unbounded simple curves, then $D$ contains at least one logarithmic tract and at most $m-1$ critical points according to multiplicity.
\end{theorem}

We now recall the classification of singularities of the inverse function due to Iversen \cite{Iversen}, using terminology found in \cite{Tracts}. Let $f$ be an entire function and consider $a \in \widehat{\mathbb{C}}$. For $r>0$, let $U_r$ be a component of $f^{-1}(D(a,r))$ (where $D(a,r)$ is the
open disc centered at $a$ with radius $r$ with respect to the spherical metric) chosen so that $r_1< r_2$ implies that $U_{r_1} \subset U_{r_2}$. Then either $\bigcap_{r} U_r = \{z\}$ for some unique $z \in \mathbb{C}$ or $\bigcap_{r} U_r = \emptyset$. These sets $U_r$ are called \textit{tracts} for $f$. 

In the first case, we have that $a=f(z)$ and $a$ is an \textit{ordinary point} if $f'(z) \neq 0$, or $a$ is a \textit{critical value} if $f'(z)=0$. In the second case,  we have a \textit{transcendental singularity} over $a$. The transcendental singularity is called \textit{direct} if there exists $r>0$ such that $f(z) \neq a$ for all $z \in U_r$. Otherwise it is \textit{indirect}. Further, a direct singularity over $a$ is called \textit{logarithmic} if $f: U_r \rightarrow D(a,r) \setminus \{a\}$ is a universal covering.

There are many functions with direct tracts which, while not logarithmic, do in fact contain a logarithmic tract, as in \Cref{InftyCurves}. This containment corresponds to an access to a logarithmic singularity in a direct tract. An example of this is given in \Cref{LogInside}, as well as in an example of Bergweiler and Eremenko \cite{BE08} of a function with no logarithmic singularities over any finite value; see \Cref{BEExample}. Further, using \Cref{Single Curve} we show that this function is in fact in the class $\mathcal{B}$. This gives an example of a function in the class $\mathcal{B}$ with no logarithmic singularities over any finite value, which was not previously known.

The organization of this paper is the following. Section 2 will be devoted to the proofs of  \Cref{Single Curve} and \Cref{InftyCurves}. Section 3 contains three examples to illustrate applications of \Cref{Single Curve} and \Cref{InftyCurves}.

\section{Proofs of \Cref{Single Curve} and \Cref{InftyCurves} }
This section is devoted to the proofs of \Cref{Single Curve} and \Cref{InftyCurves}, in which we build on some of the ideas in the proof of \cite[Theorem 5]{BE08}. First, we give a proof of \Cref{Single Curve}, which states that tracts bounded by a simple curve are logarithmic. 

\begin{proof}[Proof of \Cref{Single Curve}]
Let $\phi : \mathbb{D} \rightarrow D$ be a Riemann map, where  $\mathbb{D}$ denotes the open unit disc. The following construction is illustrated in \Cref{Phi}. The set $D$ is a Jordan domain in the Riemann sphere with boundary $\partial D \cup \{\infty\}$, so $\phi$ extends continuously and one-to-one to $\partial \mathbb{D}$, by Carath\'eodory's Theorem (\cite{CL66} and \cite[Theorem 3.1]{Harmonic}), and without loss of generality $\phi(1) = \infty$. 
So,
\[u(t)=\log \frac{|f(\phi(t))|}{R}, ~\text{for}~  t \in  \mathbb{D},\]
is a positive harmonic function in $\mathbb{D}$ with $u(t) = 0$, $t \in \partial \mathbb{D} \setminus \{1\}$. Therefore $u$ is a positive multiple of the Poisson kernel in $\mathbb{D}$ with singularity at $1$. Hence,
\[u(t) = c \operatorname{Re}\left(\frac{1+t}{1-t}\right), ~ \text{where $c>0$.}\]

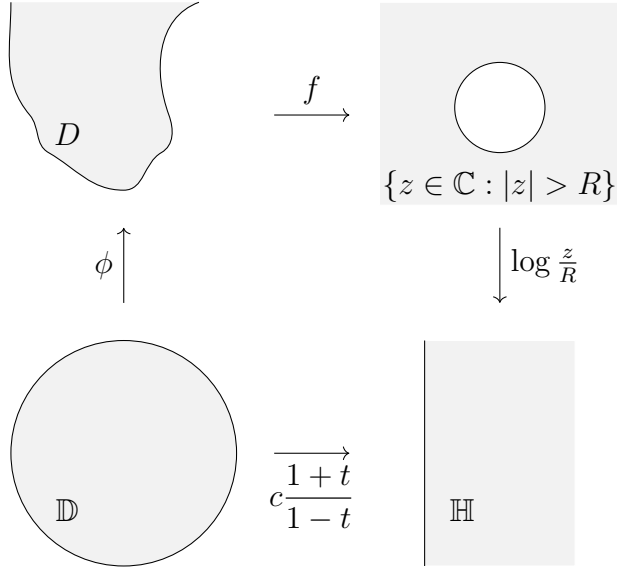
\begin{figure}
\begin{center} 
\begin{tikzpicture}[scale=1.0] 
\filldraw[ fill=black!5] (-1.5,6) to[out=-90,in=130] (-1.25,4.5) to[out=-50,in=150]  (-1,4.0) to [out=-30,in=180]  (0,3.5) to[out=0,in=-150]  (0.5,4.0) to[out=30,in=-70]  (0.6,4.5) to[out=110,in=200]  (1,6);
\draw[fill=black!5] (0,0) circle (1.5cm);
\draw[->] (0,2) -- (0,3);
\draw[->] (2,4.5) -- (3,4.5);
\draw[->] (5,3) -- (5,2);
\draw[->] (2,0) -- (3,0);

\draw[fill=black!5, draw=white]  (6.6,6) rectangle (3.4,3.3);
\draw[fill=white] (5,4.6) circle (.6 cm);

\draw[fill=black!5, draw= white]  (4.0,1.5) rectangle (6.0,-1.5);
\draw(4.0,1.5)--(4.0,-1.5);

\node at (-0.75,4.25) {${D}$};
\node at (5,3.55) {$ \{z \in \mathbb{C}: |z|> R\}$};
\node[above] at (2.5,4.5) {${f}$};
\node[left] at (0,2.5) {${\phi}$};
\node[right] at (5,2.5) {${\log \frac{z}{R}}$};
\node at (-0.75,-0.75) {$\mathbb{D}$};
\node at (4.5,-0.75) {$\mathbb{H}$};
\node[below] at (2.5,0) {$\displaystyle c \frac{1+t}{1-t}$};
\label{Commute}
\end{tikzpicture}
\caption{Construction in the proof of \Cref{Single Curve}}
\label{Phi}
\end{center}
\end{figure}

Now we can define an analytic branch $g$ of $\log f/R$ in $D$ by the monodromy theorem, since $D$ is simply connected and any local branch $g$ of $\log f/R$ can be analytically continued along any path in $D$. Then,
\[g(\phi(t))=\log  \frac{f(\phi(t))}{R}, ~\text{for}~t \in \mathbb{D},\]
is analytic in $\mathbb{D}$, with \[\operatorname{Re} g(\phi(t))=\log \frac{|f(\phi(t))|}{R}=c \operatorname{Re} \left(\frac{1+t}{1-t}\right).\] 
Hence, for some constant $\theta \in \mathbb{R}$, we have
\begin{align}
g(\phi(t))=c\frac{1+t}{1-t} + i \theta &\implies f(\phi(t))=Re^{i\theta} \exp\left(c \frac{1+t}{1-t}\right), ~\text{for}~t \in \mathbb{D}, \notag\\
&\implies f(z)=Re^{i\theta} \exp\left(c \frac{1+\phi^{-1}(z)}{1-\phi^{-1}(z)}\right), ~ \text{for}~ z \in D. \label{f}
\end{align}

It follows immediately that $f$ has no critical points in $D$. Also, there are no asymptotic paths in $D$ with finite asymptotic values, as the exponential function has none in $\mathbb{H}=\{z:\operatorname{Re}z>0\}$. Indeed, if $\gamma \rightarrow \infty$ in $D$, then $\phi^{-1}(z) \rightarrow 1$ along $ \gamma$. So, $ (1+\phi^{-1}(z))/(1-\phi^{-1}(z)) \rightarrow \infty$ in $\mathbb{H}$ as $z \rightarrow \infty$ for $z \in \gamma$. Hence, by \eqref{f}, $f(z)$ cannot tend to a finite limit as $z \rightarrow \infty$ for $z \in \gamma$. 
\end{proof}

\begin{proof}[Proof of \Cref{InftyCurves}]
Consider a direct tract $D$ bounded by more than one unbounded curve with $|f(z)|=R$, for $z\in \partial D$. If it exists, choose some channel, $G$ of $D$. Then either $|f(z)|$ will be unbounded or bounded within this channel. Let $\phi: \mathbb{D} \rightarrow G$ be a Riemann map with $\phi(z)\rightarrow \infty$ as $z \rightarrow 1$ for $z \in \mathbb{D}$. This is possible by the definition of a channel. Then, $\phi$ extends continuously and one-to-one to $\partial \mathbb{D}$ by Carath\'eodory's Theorem, and once again $\phi(1)= \infty$.

 Let $E$ be the subset of $\partial G$ where $|f(z)| \neq R$. Then,
\[u(t)=\log \frac{|f(\phi(t))|}{R}, ~\text{for}~ t \in \mathbb{D},\]
is a positive harmonic function in $\mathbb{D}$. Also, the set $\phi^{-1}(E)$ is contained in a closed arc of $\partial \mathbb{D}$ which does not contain $1$.

 First, assume $f$ is unbounded in $G$ and denote by $P(t, \zeta)$ the Poisson kernel in $\mathbb{D}$ with singularity at $\zeta$. Since $f$ is unbounded in $G$, it follows that
\begin{equation}\label{Poisson1}
u(t)=cP(t,1) + \int\displaylimits_{\phi^{-1}(E)} P(t,\zeta)u(\zeta) d\lambda(\zeta), ~\text{for}~ t \in \mathbb{D},
\end{equation}
where $c>0$ and $\lambda(\zeta)$ is the normalized Lebesgue measure on $\partial \mathbb{D}$.
Now, choose $\epsilon >0$ and $\delta >0$ such that 
\begin{equation}\label{Poisson2}
\int\displaylimits_{\phi^{-1}(E)} P(t,\zeta)u(\zeta) d\lambda(\zeta) < \epsilon,
\end{equation}
for $|t-1|<\delta$ and $t \in \mathbb{D}$. Then, choose $R_2> R_1 > R$ such that the horodiscs \[H_j=\{t\in \mathbb{D}: cP(t,1)>R_j\}, ~\text{for}~ j=1,2,\] lie inside $D(1, \delta)$ and $R_2 > R_1 + 2 \epsilon$. So, for $t \in H_2$, $u(t) > c P(t, 1) > R_2$. 
Let $\Omega$ be the component of $\{t \in \mathbb{D}: u(t) > R_2\}$ that contains $H_2$. For $t \in \Omega$, $u(t) > R_2$, so 
\[c P(t, 1) > R_2 -\int\displaylimits_{\phi^{-1}(E)} P(t,\zeta)u(\zeta) d\lambda(\zeta).\]
 Hence, by \eqref{Poisson1} and \eqref{Poisson2}, $c P(t, 1) > R_2-\epsilon > R_1$, for $t\in \Omega$. Therefore, $H_2 \subset \Omega \subset H_1$, as shown in \Cref{Horocycle}.
\begin{figure}
\begin{center} 
\begin{tikzpicture}[scale=1.0] 
\filldraw[ fill=black!5] (2.5,0) to[out=90,in=-40] (2.2,.7) to[out=140,in=-20] (1.5,.9)to[out=160,in=40] (1.1,.6)to[out=220,in=80] (1.1,0) to[out=260,in=80] (.8,-.2)to[out=260,in=120] (1.1,-.4) to[out=300,in=145] (1.6,-.7) to[out=325,in=160] (2.0,-.6)to[out=340,in=150] (2.2,-.5)to[out=320,in=270] (2.5,0);
\draw (0,0) circle (2.5cm);
\draw (1,0) circle (1.5cm);
\draw (2,0) circle (.5cm);
\node at (1.3,0) {$\Omega$};
\node at (0.25,0) {$H_1$};
\node at (2,0) {$H_2$};
\node at (-1.25,-1.25) {\large$\mathbb{D}$};
\end{tikzpicture}
\caption{$H_2 \subset \Omega \subset H_1 \subset \mathbb{D}$}
\label{Horocycle}
\end{center}
\end{figure}
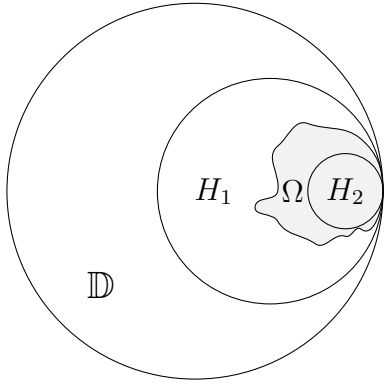

We next claim that $\Omega$ is bounded by a single curve if $R_2$ is sufficiently large. To prove this, consider the Riemann map $\psi: \mathbb{D} \rightarrow \mathbb{H}$ given by $\psi(t)=\frac{1+t}{1-t}$. Then, consider the positive harmonic function $U$ on $\mathbb{H}$ defined by
 \[U(x+iy)=u(\psi^{-1}(x+iy))=cx +x\int\displaylimits_{\psi(\phi^{-1}(E))} \frac{ds}{x^2+(y-s)^2}, ~\text{for}~ x+iy \in \mathbb{H}. \]
We have
\begin{align*}
\frac{\partial U(x+iy)}{\partial x}& = c +\int\displaylimits_{\psi(\phi^{-1}(E))} \frac{ds}{x^2+(y-s)^2}+x\frac{\partial}{\partial x}\int\displaylimits_{\psi(\phi^{-1}(E))} \frac{ds}{x^2+(y-s)^2} \\
&\geq c -\int\displaylimits_{\psi(\phi^{-1}(E))} \frac{2x^2}{(x^2+(y-s)^2)^2}ds\\
&\geq c -\frac{2}{x^2}\int\displaylimits_{\psi(\phi^{-1}(E))}ds\\
&>0
\end{align*}
for $x$ sufficiently large and $y \in \mathbb{R}$, since $\psi(\phi^{-1}(E))$ is contained in a bounded interval of the imaginary axis. So, ${\frac{\partial U(x+iy)}{\partial x}>0}$ in the half-plane $ \psi(H_1)$ for sufficiently small $\delta$. Therefore, $U(x+iy)$ is monotonic with respect to $x$ in the half-plane $ \psi(H_1)$ for any fixed~$y$. Hence $\{x+iy:U(x+iy) > R_2 \}$ is bounded by a single simple unbounded curve, and so ${\Omega = \psi^{-1}(\{x+iy:U(x+iy) > R_2 \})}$ is a Jordan domain bounded by a simple curve. Thus, by \Cref{Single Curve}, $\phi(\Omega)$ is a logarithmic tract in $G$.

Now, assume that $f$ is bounded in $G$. Then there exists $K$ such that $u(t) < K$ for $t$ in a neighborhood of the boundary singularity at $1$. On the boundary of $\mathbb{D}$, $u(t)\equiv 0$ for $t$ in a neighborhood of $1$ except possibly at $1$.  So, by the extended maximum principle \cite[Theorem 5.16]{SF1}, $u$ has boundary value $0$ at $1$. We want to show that $f(\phi(t))\rightarrow \alpha$ as $t \rightarrow 1$, where $|\alpha|=R$. Since $u\equiv0$ on $\partial \mathbb{D}$ in a neighborhood of $1$, we deduce, by the reflection principle \cite[Example 1, p. 35]{SF1}, that we can find a neighborhood, $N$ say, of $1$ in $\mathbb{C}$ to which $u$ extends harmonically. Therefore, there exists a complex conjugate $v$ of $u$ so that $u+iv$ is analytic on this neighborhood. Let 
\[g(z)= R \exp(u(\phi^{-1}(z))+iv(\phi^{-1}(z))),~\text{for}~ z \in \phi(N \cap \mathbb{D}).\]
 Then, $g$ is analytic on $\phi(N \cap \mathbb{D}) \subset G$ and $|g(z)|= |f(z)|$ in $\phi(N \cap \mathbb{D})$. Non-constant analytic maps are open, so $g(z)= c f(z)$ in $\phi(N \cap \mathbb{D})$, where $|c|=1$. Therefore, $\arg f(z)$ tends to a finite limit as $z \rightarrow \infty$ in $\phi(N \cap \mathbb{D})$, and hence in $G$, since any sequence tending to infinity in $G$ is eventually contained in $\phi(N \cap \mathbb{D})$. Therefore, there exists $\alpha$ such that $f(z) \rightarrow \alpha$ as $z \rightarrow \infty$ in $G$, with $|\alpha|=R$. 
 
 Finally, if $D$ is bounded by $m$ distinct unbounded simple curves, then following the method in the proof of \Cref{Single Curve}, we again have
\[u(t)=\log \frac{|f(\phi(t))|}{R}, ~\text{for}~  t \in  \mathbb{D},\]
where $\phi: \mathbb{D} \rightarrow D$ is a conformal map and each access to infinity in $D$ corresponds under $\phi$ to a family of paths tending to a point on $ \partial \mathbb{D}$ along which $u$ tends to an asymptotic value, either $0$ or $\infty$. Suppose there exist $n \leq m$ accesses to infinity on which $f$ is unbounded. Then, there exist $n$ points $ \zeta_1, \ldots, \zeta_n \in \partial \mathbb{D}$ such that $u$ is a positive harmonic function in $\mathbb{D}$ with $u(t) = 0$, for $t \in \partial \mathbb{D} \setminus \{ \zeta_1, \ldots, \zeta_n\}$. Therefore $u$ is a sum of positive multiples of the Poisson kernel in $\mathbb{D}$ with singularities at $ \zeta_1, \ldots, \zeta_n$. Hence,
\[u(t) = \sum_{k=1}^n c_k \operatorname{Re}\left(\frac{\zeta_k+t}{\zeta_k-t}\right), ~ \text{where $c_k>0$, for $k=1, \ldots, n$.}\]
Now, again, we can define an analytic branch $g$ of $\log f/R$ in $D$ by the monodromy theorem, since $D$ is simply connected and any local branch $g$ of $\log f/R$ can be analytically continued along any path in $D$. Then,
\[g(\phi(t))=\log  \frac{f(\phi(t))}{R}, ~\text{for}~t \in \mathbb{D},\]
is analytic in $\mathbb{D}$, with \[\operatorname{Re} g(\phi(t))=\log \frac{|f(\phi(t))|}{R}=\sum_{k=1}^n c_k \operatorname{Re}\left(\frac{\zeta_k+t}{\zeta_k-t}\right).\] 
Hence, for some constant $\theta \in \mathbb{R}$, we have
\begin{align}
g(\phi(t))=\sum_{k=1}^n c_k \frac{\zeta_k+t}{\zeta_k-t} + i \theta &\implies f(\phi(t))=Re^{i\theta} \exp\left(\sum_{k=1}^n c_k \frac{\zeta_k+t}{\zeta_k-t}\right), ~\text{for}~t \in \mathbb{D}. \notag
\end{align}
The critical points of $f(\phi(t))$ for $t \in \mathbb{D}$ are the solutions of 
\[ \sum_{k=1}^n c_k \frac{2\zeta_k}{(\zeta_k-t)^2}=0,\]
for which there are at most $2n-2$ solutions in $\mathbb{C}$. Further, by the Cauchy-Riemann equations, critical points of $f(\phi(t))$ occur if and only if $|\nabla u(t)|=0$, where $\nabla u$ denotes the gradient of $u$. For a point $t \in \mathbb{D}$, noting that $\zeta_k= 1/ \overline{\zeta}_k$ gives that
\[ \operatorname{Re}\left(\frac{\zeta_k+t}{\zeta_k-t}\right) =  - \operatorname{Re}\left(\frac{\zeta_k+1/\overline{t}}{\zeta_k-1/\overline{t}}\right), ~\text{for}~ k = 1, \ldots, n.\]
  Hence, each Poisson kernel is `symmetric' under reflection in the unit circle and so a sum of Poisson kernels is as well. Further, no solutions of $|\nabla u(t)|=0$ lie on $\partial \mathbb{D}$, by the behavior of the Poisson kernel near $\partial \mathbb{D}$. So, the solutions of $|\nabla u(t)|=0$ occur in pairs and are symmetric with respect to $\partial\mathbb{D}$. Therefore, $\mathbb{D}$ contains at most $n-1$ critical points of $f(\phi(t))$ and hence $D$  contains at most $n-1$ critical points of $f$.
\end{proof}

\section{Examples}
In this section, we give three examples to show the kinds of direct tracts that can exist, and to which we can apply \Cref{Single Curve} and \Cref{InftyCurves}. In our first example, we also show that our results give a function in class $\mathcal{B}$ which has no logarithmic singularities over any finite value. 

\begin{figure}[!hbt]
  \centering
    \fbox{\includegraphics[width=0.4\textwidth]{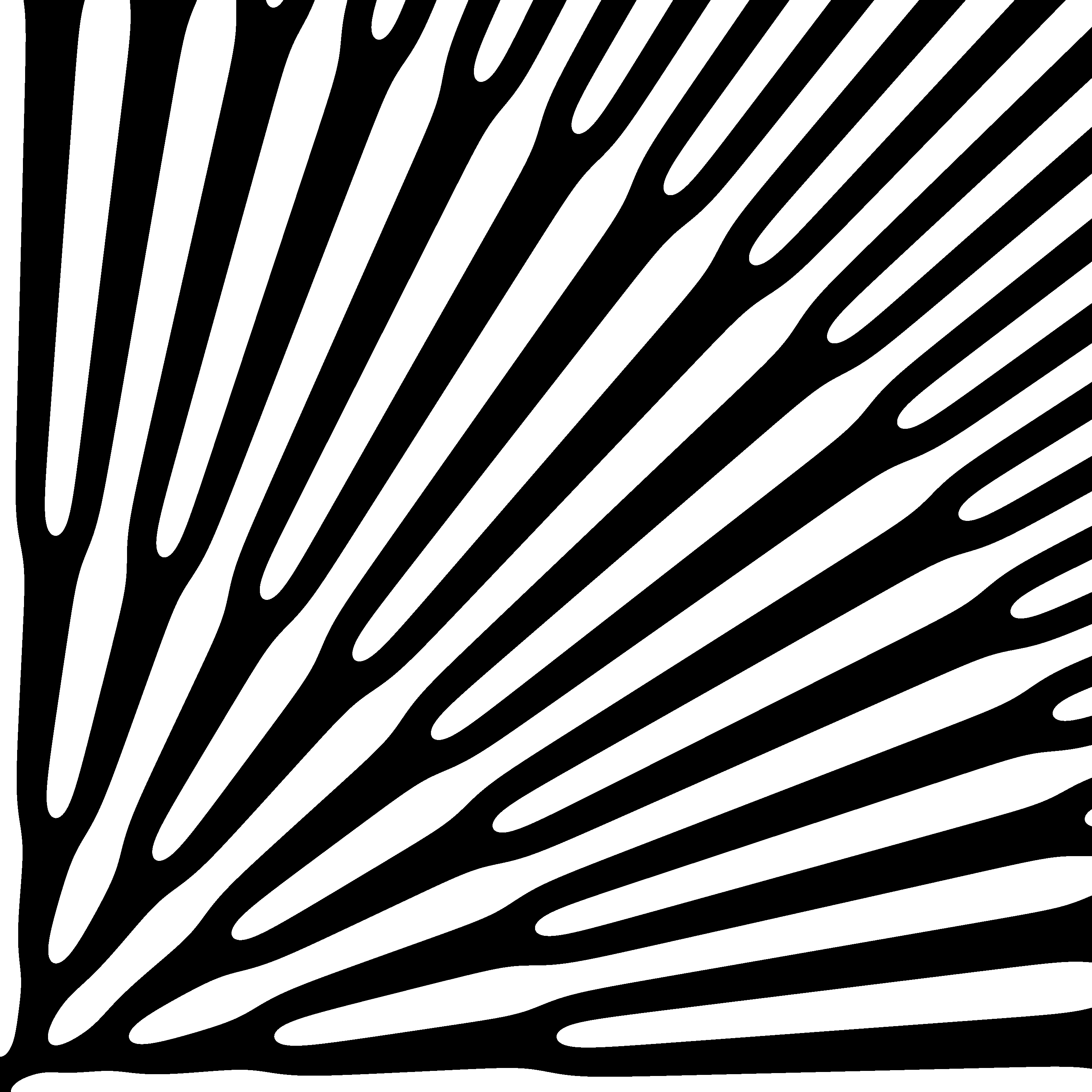}}
  \caption{The direct tracts of $h$ in white for $0\leq\operatorname{Re} z \leq 250$ and $0\leq\operatorname{Im} z\leq 250$.}
\end{figure}
\begin{example}\label{BEExample}
Consider the entire function
\[h(z)=\exp(g(z)), ~\text{where} ~g(z)=\sum_{k=1}^{\infty} \left(\frac{z}{2^k}\right)^{2^k}.\]
Then, $h$ has infinitely many direct singularities, but no logarithmic singularity over any finite value, and is in the class $\mathcal{B}$.

The first two statements are proved in \cite{BE08}, so it remains to check that $h$ is in the class $\mathcal{B}$. In order to see this, we first observe that the only possible finite asymptotic value of $h$ is $0$, since the only finite asymptotic value of the exponential function is $0$. We will show that every direct tract (over $\infty$) for some fixed boundary value is bounded by a single curve, and so all the direct tracts are logarithmic by \Cref{Single Curve}. This will imply that no critical points lie in these tracts and so the critical values of $h$ form a bounded set. We can then conclude that $h$ is in the class $\mathcal{B}$.

Following the notation and construction in \cite{BE08}, fix $\varepsilon$ with $0<\varepsilon\leq \frac{1}{8}$ and set $r_n=(1+\varepsilon)2^{n+1}$ and $r_n'=(1-2\varepsilon)2^{n+2}$ for $n\in \mathbb{N}$. Then for $j\in\{0,1,\cdots,2^n-1\}$ we define the sets

\[B_{j,n}=\left\{r\exp\left(\frac{\pi i}{2^n}+\frac{2\pi ij}{2^n}\right) :r_n\leq r\leq r'_n \right\}\]
and
\[C^{\pm}_{j,n}=\left\{r \exp\left(\frac{\pi i}{2^n}+\frac{2 \pi i j}{2^n} \pm \frac{r-r'_n}{r_{n+1}-r'_n} \frac{\pi i}{2^{n+1}} \right) : r'_n\leq r\leq r_{n+1} \right\}.\]

Bergweiler and Eremenko \cite[Section 6]{BE08} show that every unbounded simple path starting at $0$ and lying in the infinite tree,
\[T=[-ir_1,ir_1] \cup \bigcup^\infty_{n=1} \bigcup^{2^n-1}_{j=1}(B_{j,n} \cup C^{\pm}_{j,n} ), \] is an asymptotic curve along which the function tends to 0. Further,
\[\operatorname{Re} g(z) < -2^{2^n} ~\text{for}~ z \in B_{j,n} \cup C^{\pm}_{j,n}.\]
So, if some direct tract (over $\infty$) with a sufficiently large boundary value was not bounded by a single curve, then at least one of the following three possibilities would occur:
\begin{enumerate}
\item there would be another direct tract over 0,
 \item the tract would have to cross the infinite tree $T$,
 \item there would be zeros of the function $h$. 
  \end{enumerate}
  The first case is shown not to happen in \cite[Section 6]{BE08} by proving that $\arg g(r e^{i \theta})$ is an increasing function of $\theta$, that it increases by $2^n 2 \pi$ as $\theta$ increases by $2 \pi$, and then using a counting argument to show that all the direct tracts are accounted for and there can be no others. We can assume that the tract does not cross the infinite tree $T$, since $h(z)$ is bounded on $T$ by a value smaller than the tract boundary value. Finally, the last case cannot happen since the exponential function omits $0$.  Therefore, each direct tract (over $\infty$) is bounded by a single curve, so by \Cref{Single Curve} we can deduce that $h\in \mathcal{B}$.
\end{example}

Next, we illustrate \Cref{InftyCurves} by giving an example of a direct non-logarithmic tract with infinitely many logarithmic tracts inside it.

\begin{figure}[!h]
  \centering
    \fbox{\includegraphics[width=0.4\textwidth]{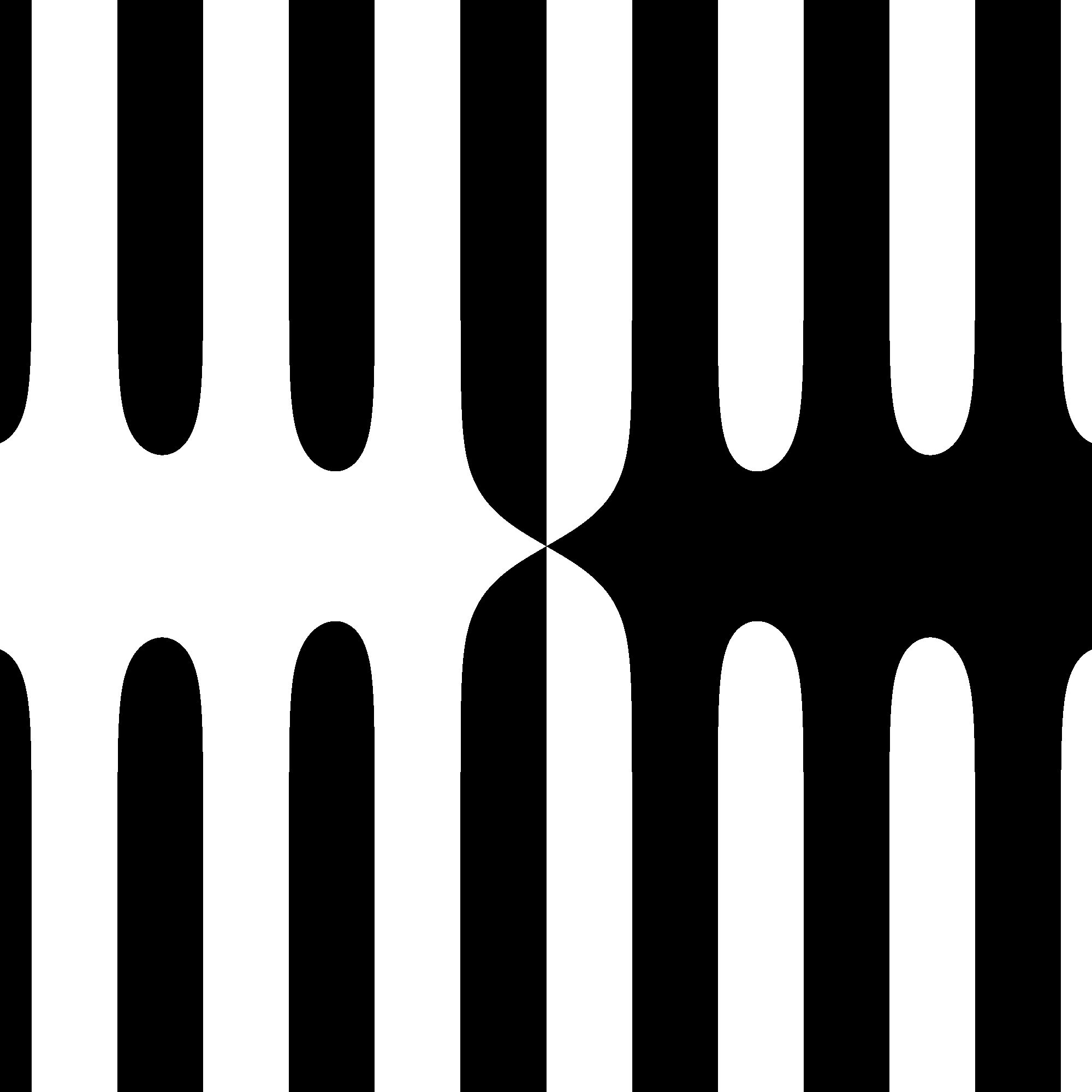}}
  \caption{The tracts of $\exp(\sin(z)-z)$ in white.}
  \label{fig:LogInside}
\end{figure}
\begin{example} \label{LogInside}
Consider $f(z)=\exp(\sin(z)-z)$. First, $f(z) \rightarrow \infty$ as $z \rightarrow \infty$ along the negative real axis and $f(z) \rightarrow \infty$ as $z \rightarrow \infty$ along translates of the imaginary axis by $\frac{\pi}{2}+2k\pi$, for $k \in \mathbb{Z}$. Also, $f(z) \rightarrow 0$  as $z \rightarrow \infty$ along the positive real axis and along translates of the imaginary axis by $\frac{3\pi}{2}+2k\pi$, for $k  \in \mathbb{Z}$. So, $f$ will have infinitely many direct tracts in the right half-plane and one direct tract in the left half-plane, with $|f(z)|=R$ on $\partial D$ for some suitable $R>0$. See \Cref{fig:LogInside}, where $R=1$. Further, $f$ has no zeros and its only finite asymptotic value is $0$,  along the positive real axis and translates of the imaginary axis by $\frac{3\pi}{2}+2k\pi$, for $k \in \mathbb{Z}$.  Hence, by \Cref{InftyCurves}, the direct tract in the left half-plane contains infinitely many logarithmic tracts, each corresponding to channels about  translates of the positive and negative imaginary axes by $\frac{\pi}{2}+2k\pi$, for $k \in \mathbb{Z}$. In contrast, the tracts in the right half-plane are all logarithmic tracts.
\end{example}

Finally, to illustrate  \Cref{InftyCurves}, we give a simple example of a transcendental entire function with a simply connected direct tract bounded by finitely many boundary curves.
\begin{example} \label{Finite}
\begin{figure}[!hbt]
  \centering
    \fbox{\includegraphics[width=0.4\textwidth]{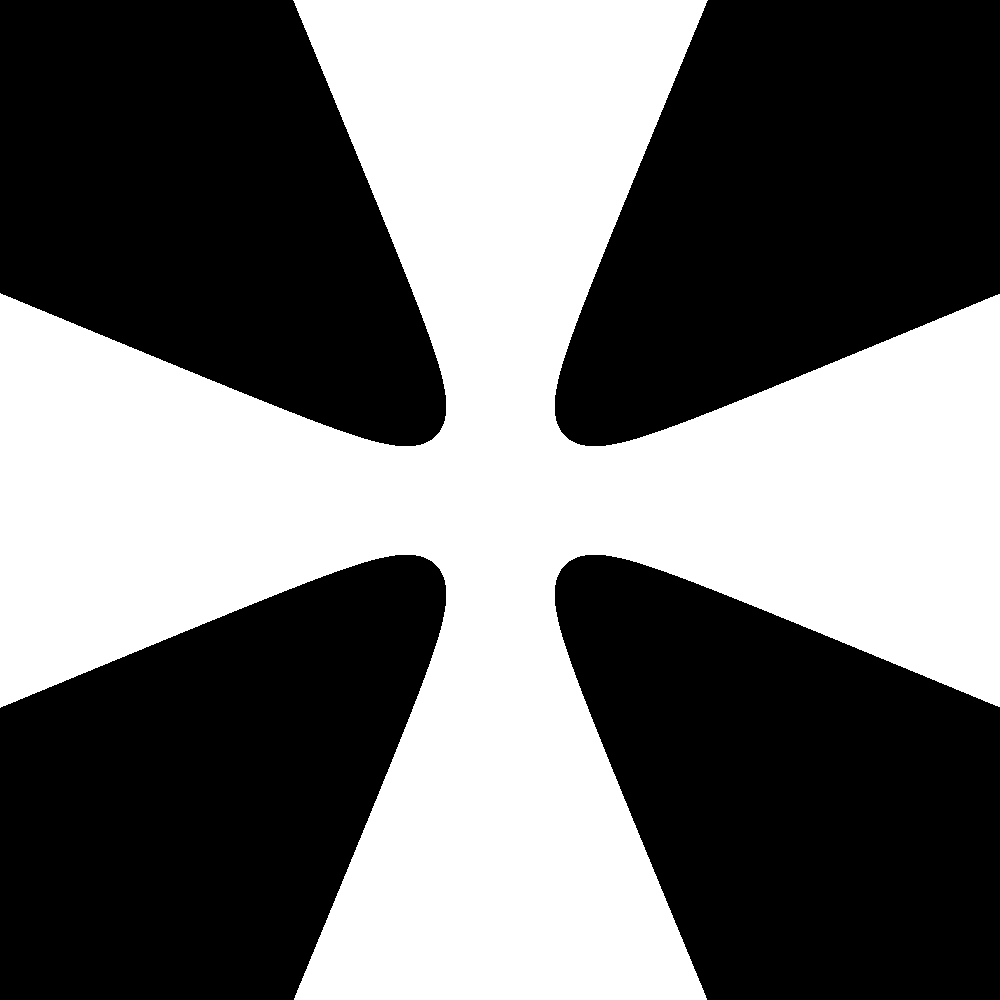}}
  \caption{The tract of $2\exp(z^4)$ in white.}
  \label{exp(z4)}
\end{figure}
Consider $G(z)=2\exp(z^4)$. 
First, $G(z) \rightarrow \infty$ as $z \rightarrow \infty$ along the real and imaginary axes, and $G(z) \rightarrow 0$ as $z \rightarrow \infty$ along the rays with angle an odd multiple of $\pi/4$. With a boundary value of $R=1$, as in \Cref{exp(z4)}, a neighborhood about the origin is contained in the direct tract. Further,  on the lines with angle an odd multiple of $\pi/8$ the modulus of $G$ is~$2$. Hence, $G$ has one direct tract bounded by four unbounded simple curves. By \Cref{InftyCurves}, $G$ contains at most three critical points, and in fact,  contains a single critical point of multiplicity~$3$ at $0$. The direct tract of $G$ has four channels on which $f$ is unbounded and hence contains four logarithmic tracts.
\end{example}
\subsection*{Acknowledgements}
I would like to thank my supervisors Prof. Phil Rippon and Prof. Gwyneth Stallard for their guidance and encouragement in the preparation of this paper. I would also like to thank Dr. Dave Sixsmith for mentioning the problem of finding an example of a function in the class $\mathcal{B}$ with a direct non-logarithmic singularity over a finite value, which led to \Cref{BEExample}.
\bibliographystyle{abbrv}
\bibliography{mybib}{}
\end{document}